\documentclass[12pt]{article}

\usepackage{amssymb}
\usepackage{amsthm}
\usepackage{amsmath}

\newcommand{\EE}{\mathbb{E}}

\newcommand{\RR}{\mathbb{R}}

\newtheorem{theorem}{Theorem}

\newcommand{\ep}{\epsilon}

\newcommand{\lbr}{\langle}
\newcommand{\rbr}{\rangle}

\newcommand{\dsp}{\displaystyle}

\author{Jean-Louis Krivine\\
\footnotesize{Universit\'e Paris-Cité, C.N.R.S.}\\
\footnotesize{krivine@irif.fr}}

\title{A note about Grothendieck's constant}
\date{\normalsize May 24, 2023}

\begin{document}\noindent
\maketitle\noindent
\emph{Grothendieck's constant}, denoted $K_G$, is defined as the smallest real $K>0$ such that we can
write $\lbr\vec{x},\vec{y}\rbr=K\,\EE(\mbox{sign}(U_{\vec{x}})\mbox{sign}(V_{\vec{y}}))$ with
$\vec{x},\vec{y}\in\mathbb{S}$ (the unit sphere of $l^2$) and $U_{\vec{x}},V_{\vec{y}}$ are random variables,
which are measurable functions of $\vec{x},\vec{y}$ respectively.\\
Its existence was proved by A. Grothendieck in \cite{groth} where it is shown that~:\\
$\pi/2\le K_G\le\mbox{sh}(\pi/2)=2.301\ldots$\\
This constant is important in various areas such as functional analysis, algorithmic complexity and
quantum mechanics~: see \cite{pisier}. Its exact value is unknown.\\
In \cite{kri}, it is shown that $K_G\le\pi/2\ln(1+\sqrt{2})=1.782\ldots$ This result is improved in the prominent article \cite{nao} which proves~:
\begin{equation}\label{strict}
K_G<\frac{\pi}{2\ln(1+\sqrt{2})}
\end{equation}
Of course, the method gives a better upper bound but the authors did not considered it useful
to give it explicitly for the moment. The important fact in~(\ref{strict}) is the symbol $<$.\\[3pt]
As explained below, the proof given in \cite{nao} is divided into two parts and the aim of the present note
is to give another proof of the first part.\\[3pt]
Let $(X_i,Y_i)(0\le i\le n-1)$ be independent pairs of centered gaussian normal random variables,
such that $\EE(X_iY_i)=t$.\\
Let $\vec{X}=(X_0,\ldots,X_{n-1})$, $\vec{Y}=(Y_0,\ldots,Y_{n-1}),
\vec{x}=(x_0,\ldots,x_{n-1})$, $\vec{y}=(y_0,\ldots,y_{n-1})$.\\
Let $F,G:\RR^n\to\RR$ be two odd (i.e. $F(-\vec{x})=-F(\vec{x})$) measurable functions. We set~:\\
$$\dsp\Phi_{F,G}(t)=\EE\left[\mbox{sign}(F(\vec{X}))\mbox{sign}(G(\vec{Y}))\right]$$
or else~:
\begin{equation}\label{efxgyt}
\Phi_{F,G}(t)=(2\pi\sqrt{1-t^2})^{-n}\int_{\mathbb{R}^{2n}}
\mbox{sign}(F(\vec{x}))\mbox{sign}(G(\vec{y}))
e^{-\frac{\|\vec{x}\|^2+\|\vec{y}\|^2-2t\lbr\vec{x},\vec{y}\rbr}{2(1-t^2)}}d\vec{x}d\vec{y}
\end{equation}
$\Phi_{F,G}(t)$ is an odd function of $t$, which is analytic around $0$, in fact for $|\mbox{Re}(t)|<1$.\\
Therefore $\Phi_{F,G}(i)/i$ is real. We have~:
\begin{equation}\label{iphi}
\Phi_{F,G}(i)/i=(2\pi\sqrt{2})^{-n}\int_{\mathbb{R}^{2n}}
\mbox{sign}(F(\vec{x}))\mbox{sign}(G(\vec{y}))e^{-\frac{\|\vec{x}\|^2+\|\vec{y}\|^2}{4}}
\sin\frac{\lbr\vec{x},\vec{y}\rbr}{2}d\vec{x}d\vec{y}
\end{equation}
The first part of the proof of (\ref{strict}), which is section~4 of~\cite{nao},
consists in showing the~:
\begin{theorem}\label{contr}\ \\
There exists an integer $n\ge1$ and two odd functions $F,G:\RR^n\to\RR$ such that~:\\
\centerline{$\dsp\Phi_{F,G}(i)/i>\frac{2}{\pi}\ln(1+\sqrt{2})$.}
\end{theorem}\noindent
H.~K\"onig has shown that $n$ must be $>1$~; a proof of this is given in section~6 of \cite{nao}.\\[3pt]
In \cite{nao}, the authors choose $n=2$~; $F(\vec{X})=G(\vec{X})=X_0+\ep H_5(X_1)$ where $H_n(x)$ is
the Hermite polynomial of degree $n$~; we have $H_5(x)=x^5-10x^3+15x$~;
$\ep$ is a positive real which decreases to $0$.\\[3pt]
Here we take $n=3$ with~:\\
$F(\vec{X})=X_1\cos(\ep H_2(X_0))+X_2\sin(\ep H_2(X_0))$ and
$G(\vec{Y})=Y_1\cos(\ep H_2(Y_0))-Y_2\sin(\ep H_2(Y_0))$.\\
where $\ep$ is a \emph{fixed} positive real and $H_2(x)=x^2-1$.\\
Applying the well known formula \ $\dsp \EE(\mbox{sign}(X)\mbox{sign}(Y))=
\frac{2}{\pi}\mbox{Arcsin}(\EE(XY))$, we get~:\\
$$\EE\left[\mbox{sign}(F(\vec{X}))\mbox{sign}(G(\vec{Y}))\right]
=\frac{2}{\pi}\EE\left[\mbox{Arcsin}\left(t\,\cos(\ep(X_0^2+Y_0^2-2))\right)\right]$$
or else~:
$$\Phi_{F,G}(t)=\frac{2}{\pi}\int_{\mathbb{R}^{2}}\mbox{Arcsin}\left(t\cos(\ep(x^2+y^2-2))\right)
e^{-\frac{x^2+y^2-2txy}{2(1-t^2)}}\frac{dx\,dy}{2\pi\sqrt{1-t^2}}$$
Let $t=i$ and $\eta=2\ep$~:
$$\Phi_{F,G}(i)/i=\frac{2}{\pi}\int_{\mathbb{R}^{2}}\mbox{Argsh}
\left(\cos(\ep(x^2+y^2-2))\right)e^{-\frac{x^2+y^2}{4}}\cos\left(\frac{xy}{2}\right)
\frac{dx\,dy}{2\pi\sqrt{2}}$$
$$=\frac{2}{\pi}\int_0^\infty\int_{-\pi}^{\pi}\mbox{Argsh}
\left(\cos(\ep(r^2-2))\right)e^{-\frac{r^2}{4}}\cos\left(\frac{r^2\sin2\theta}{4}\right)
\frac{r\,dr\,d\theta}{2\pi\sqrt{2}}$$
$$=\frac{2}{\pi}\int_0^\infty\int_{0}^{\pi}\mbox{Argsh}
\left(\cos(\eta(2\rho-1))\right)e^{-\rho}\cos\left(\rho\sin\theta\right)
\frac{\sqrt{2}}{\pi}d\rho\,d\theta$$
This integral is not difficult to compute with a suitable software, which also gives a good
value for $\eta$~;
with $\eta=0.228$ we find $\dsp 0.56161447>\frac{2}{\pi}\ln(1+\sqrt{2})=0.56109985\ldots$
\hspace{\fill}$\square$\\[3pt]
Here are the details of the calculations in Mathematica and Maxima~:

\smallskip\noindent
\emph{Computation in Mathematica}\\
\verb?e = 0.228; 2*(Sqrt[2]/Pi^2)*NIntegrate[ArcSinh[Cos[e*(2r-1)]]?\\
\verb?*Exp[-r]*Cos[r*Sin[t]],{r,0,Infinity},{t,0,Pi}]?\\
\verb?0.561614475916681?

\smallskip\noindent
\emph{Computations in Maxima}\\
\verb?e:0.228$ float((2*sqrt(2)/(%pi)^2)*romberg(romberg((%e^(-r)*cos(r*sin(t))?\\
\verb?*asinh(cos(e*(2*r-1)))),r,0,100),t,0,%pi));?\\
\verb?0.5616147985832746?\\
\verb?e:0.228$ float(2*sqrt(2)/(%pi)^2)*quad_qag(romberg((%e^(-r)*cos(r*sin(t))?\\
\verb?*asinh(cos(e*(2*r-1)))),r,0,100),t,0,%pi,3);?\\
\verb?[0.561614489984486,5.087584361677373*10^-10,8.883967107886724,0]?

\smallskip\noindent
If the inverse power series of $\Phi_{F,G}(t)$ had been alternating, we would have obtained in this way
an upper bound for $K_G$, that is \ $i/\Phi_{F,G}(i)<1,7806$. But it is not, as is easily checked,
using the same computation tools. The second part of the proof in \cite{nao} must therefore now be applied.
It occupies section 5 of this article and uses only the above theorem~\ref{contr}.


\begin{thebibliography}{99}
\bibitem{groth}A. Grothendieck. \emph{Résumé de la théorie métrique des produits tensoriels topologiques.} Bol. Soc. Mat. Sao Paulo, 8:1–79 (1953).
\bibitem{kon}H. K\"onig. \emph{On an extremal problem originating in questions of unconditional convergence.} In Recent progress in multivariate approximation (Witten-Bommerholz, 2000), vol. 137 of Internat. Ser. Num.
Math., p. 185–192. Birkh\"auser, Basel (2001).
\bibitem{nao}M. Braverman, K. Makarychev, Y. Makarychev, and A. Naor. \emph{The Grothendieck constant is
strictly smaller than Krivine’s bound.} 52nd Annual IEEE Symposium on Foundations of Computer Science
(FOCS) p. 453-462 (2011).
\bibitem{kri}J.-L. Krivine. \emph{Constantes de Grothendieck et fonctions de type positif sur les sphères.} Advances in Math. 31, p. 16-30 (1979).
\bibitem{pisier}G. Pisier. \emph{Grothendieck’s theorem, past and present.} Preprint available at
http://arxiv.org/abs/1101.4195 (2011).
\end{thebibliography}
\end{document}